\documentclass{amsart}
\usepackage[active]{srcltx}

\usepackage[colorlinks,plainpages,backref]{hyperref}

\usepackage[all]{xy }

\theoremstyle{definition}

\theoremstyle{remark}

\numberwithin{equation}{section}

\newcommand\z{\ensuremath{\mathbb Z}}
\newcommand\be{\begin{equation}}
\newcommand\ee{\end{equation}}
\newcommand\vez {\ensuremath{\times}}
\newcommand{\wh}[1]{\ensuremath{\widehat{#1}}}
\newcommand\wt[1]{\ensuremath{\widetilde{#1}}}
\newcommand\codim{\ensuremath{\operatorname{codim}\,}}
\newcommand\C{\ensuremath{\mathbb C}}
\newcommand\p[1]{\ensuremath{\mathbb{P}^{#1}}}
\newcommand\lar{\ensuremath{\,\longrightarrow\,}}
\newcommand\rar{\ensuremath{\dashrightarrow}}
\newcommand\ra{\ensuremath{\,\rightarrow\,}}
\newcommand\ls{\ensuremath{_\star}}
\newcommand\us{\ensuremath{^\star}}
\newcommand\bbb[1]{\ensuremath{\mathbb{#1}}}
\newcommand\cl[1]{\ensuremath{\mathcal{#1}}}
\newcommand\s[1]{\ensuremath{\cl F_{#1}}}
\newcommand\ov[1]{\ensuremath{\overline{#1}}}
\newcommand\ba[1]{\begin{array}{#1}}\newcommand\ea{\end{array}}
\newcommand\id[1]{\ensuremath{ \langle {#1} \rangle}}
\newcommand\ie{{\em i.e.\em, }}
\newcommand\x{\ensuremath{\mathbb X}}
\newcommand\G{\ensuremath{\mathbb G}}
\newcommand\E{\ensuremath{\mathbb E}}
\newcommand\st{\ensuremath{\,|\,}}
\newcommand\escala[1]{\setlength{\unitlength}{#1cm}}
\newcommand\bi{\begin{itemize}}\newcommand\ei{\end{itemize}}
\newcommand\ve{\ensuremath{^\vee}}
\newcommand\rf[1]{\hbox{\rm(\ref{#1})}}
\newcommand\ds[1]{\displaystyle{#1}}
\newcommand\na[1]{\noalign{\vskip#1pt}}

\newcommand\Pic{\operatorname{Pic}}

\newcommand\pd[1]{\ensuremath{\mathbb{\check{P}}^{#1}}}
\newcommand\Hom{{\rm Hom}} 
\newcommand\wtn{\ensuremath{\wt{NL}(W,d)}}
\newcommand\nld{\ensuremath{{NL}(W,d)}}
\newcommand\ch{\ensuremath{\operatorname{ch}}}
\newcommand\isom{\ensuremath{\cong}}
\newcommand\Tau{\mbox{\Large$\tau$}\!}

\begin{document}

\title[Surfaces containing 
\linebreak 
an elliptic quartic
  curve]{Enumeration of surfaces containing 
an elliptic quartic curve}

\author[F.\,Cukierman] {F.\,Cukierman$^\spadesuit$}
\address{Universidad de Buenos Aires
\\
Ciudad Universitaria, Pabell\'on 1
\\(1428) Buenos Aires Argentina}
\email{fcukier@dm.uba.ar}
\thanks{$^\spadesuit$ Partially supported by 
CONICET-Argentina.}

\author[A.F.\,Lopez] {A.F.\,Lopez$^\heartsuit$}
\address{Dipartimento di Matematica e Fisica, Universit\`a di Roma
Tre, Largo San Leonardo Murialdo 1, 00146, Roma, Italy}
\email{lopez@mat.uniroma3.it}
\thanks{$\heartsuit$\,Partially supported by PRIN Geometria delle
variet\`a algebriche e dei loro spazi di
moduli.}


\author[I.\,Vainsencher] {I.\,Vainsencher$^\star$}
\address{ICEX-Departamento de Matem\'atica-UFMG, Av. Ant\^onio Carlos, 6627 - Caixa Postal 702
\\CEP 31270-901 Belo Horizonte, MG, Brazil}
\email{israel@mat.ufmg.br} 

\thanks{%
$^\star$\,  Partially supported by CNPQ-Brasil.}

\date{\today}



\date{\today}

\dedicatory{Dedicated to Steve Kleiman on the occasion of his
  70th birthday}

\commby{Lev Borisov}

\begin{abstract}
A very general surface of degree at least four in \p3
contains no curves other than intersections with
surfaces. We find a formula for the degree of the locus
of surfaces in \p3 of degree at least five 
which contain
some elliptic quartic curves.  We also compute the
degree of the locus of quartic surfaces containing an
elliptic quartic curve, a case not covered by that formula.
\end{abstract}

\maketitle


\bibliographystyle{amsplain}
\section{Introduction}

The Noether-Lefschetz theorem asserts that all curves
contained in a very general surface $F$ of degree at
least four in \p3 are complete intersections.  This is
usually rephrased by saying that the Picard group is
\z.  Noether-Lefschetz theory shows that, roughly
speaking, each additional generator for Pic$F$
decreases the dimension of the locus of such $F$ in $\p
N=|\cl O_{\p3}(d)|,\,d\geq4$.

Let $W$ be a closed, irreducible subvariety of the
Hilbert scheme of 
curves in \p3 with Hilbert polynomial $p_W(t)$.
Let us denote by $NL(W,d)$ the subset of $\p N$
defined by the requirement that the surface contains
some member of $W$.

The purpose of this note is to address 
the question of determining  the degree of
$NL(W,d)$ for the family of elliptic
quartic curves in \p 3.

When $W$ is the family of lines, or conics, or twisted
cubics, formulas for \nld\ have been found in
\cite{aafi}.  There as here, we follow the strategy of
using Bott's formula as explained in
\cite{ell&stromme}.  We get a polynomial formula
(\ref{f}) valid for $d\geq5$.  We also compute the
degree ({\bf38475}) of the locus of quartic surfaces
containing an 
elliptic quartic curve. The case of quartic surfaces is
not covered by the formula essentially because the map
that forgets the curve shrinks dimensions: generically,
it contracts a pair of disjoint pencils, see \ref{maincor}.

\section{there is a polynomial formula}
Let $W$ be a closed, irreducible subvariety of the
Hilbert scheme of 
curves in \p3 with Hilbert polynomial $p_W(t)$.
Let 
\be\label{wpn}\xymatrix
@C1pc
{\,W\vez\p N \ar@{->}[rr]^{\hskip.51cmp_2}\ar@{->}[d]
_{p_1}^{\phantom{--}}\ 
& &\p N\\
W&&
  }
\ee
denote the projection maps from $W\vez\p N$.
Castelnuovo-Mumford regularity \cite{M}
shows that for all
$d>>0$, the subset $\wt{NL}(W,d)$ of pairs
$(C,F)$ in $W\vez\p N$ such that the curve $C$ is
contained in the surface $F$ is a projective
bundle over $W$ via $p_1$. We have
$$
\codim_{W\vez\p N}\wtn=p_W(d).
$$ 
For instance, if $W$ is the Grassmannian of lines in
\p3, then $p_W(d)=d+1$ and so 
$\dim \wtn=N-(d-3),\,d\geq1$. 

Let us denote by $NL(W,d)$ the
subset of $\p N$ defined by the
requirement that the surface contains some member
of $W$. In other words, with notation as in (\ref{wpn}),
$$
NL(W,d)=p_2(\wtn).
$$

We assume
henceforth that the general member of $W$ is a smooth
curve. 

\subsection{Proposition}
\label{pol}
{\em 
For fixed $W$ we have that $\deg{}NL(W,d)$ is a
polynomial in $d$ of degre $\leq 3\dim\, W$, for all
$d>>0$.}

\medskip
{\em Proof.}
Let $\wt C\subset W\vez\p 3$ be the
universal curve. Likewise, let $\wt F\subset\p
N\vez\p3$ be the universal surface of degree $d$.
Write \wh C,\,\wh F\ for their pullbacks to 
$W\vez\p N\vez\p3$.
We have the diagram of sheaves over $Y:=W\vez\p N
\vez\p3$,
$$
\xymatrix
{\cl O_Y\,\ar@{->}[dr]_\rho
\ar@{->}[r]&\cl O_Y(\wh F)\ \ar@{->}[r]
\ar@{->>}[d]& 
\cl O_{\wh F}(\wh F)\\
&\cl O_{\wh C}(\wh F)& 
}
$$
By construction, the slant arrow $\rho$ vanishes at a
point $(C,F,x)\in W\vez\p N \vez\p3$ if and only if
$x\in F\cap C$.  We have $C\subset F$ when the previous
condition holds for all $x\in C$ (point with values in
any \C-algebra).  Thus $\wtn$ is equal to the
scheme of zeros of $\rho$ along the fibers of the
projection $p_{12}:\wh C
\ra{}W\vez\p N$.
Recalling \cite[(2.1),p.\,14]{altman-kleiman}, this is
the same as the zeros of the adjoint section of the
direct image vector bundle 
$p_{12}{}\ls(\cl O_{\wh C}(\wh{F}))$.
Let 
$$\xymatrix
@C1pc
{W\vez\p 3\ \ar@{->}[rr]^{\hskip.621cm q_2}\ar@{->}[d]_{q_1}& &\p 3\\
W&&
  }
$$
denote the projection maps from $W\vez\p 3$.
Since $\cl O(\wt F)=\cl O_{\p N}(1)\otimes\cl O_{\p3}(d)$,
by projection formula 
we have to make do with a section of 
$\cl O_{\p N}(1)\otimes\cl E_d$, where
\be\label{ed}
\cl E_d=
q_1{}\ls(\cl O_{\wt C}(d)).
\ee
By Castelnuovo-Mumford and base change theory,
there is an integer $d_0$ such that
$\cl E_d$ is a vector bundle of rank $p_W(d)$
for all $d\geq d_0$\,(=regularity, 
see Remarks\,\ref{cmr}.).
In fact, it fits into the exact sequence of
vector bundles over $W$,
\be\label{dd}
\xymatrix
{
0\ \ar@{->}[r]& \ \ar@{=}[d]
q_1{}\ls(\cl I_{\wt C}\otimes\cl O_{\p3}(d))\,\ar@{->}[r]&
q_1{}\ls(\cl O_{\p3}(d))\,\ar@{->}[r] \ar@{=}[d]&
q_1{}\ls(\cl O_{\wt C}(d))\ \ar@{->}[r]& 0\\
& \ \,\,
\cl D_d\ \,\ar@{>->}[r]&
\s d\,\ar@{->>}[r]&
\,\cl E_d\ \ar@{=}[u]& 
}
\ee
where we set for short 
$$
\s d=H^0(\cl O_{\p3}(d)),
$$ 
(trivial vector bundle with fiber) the space of
polynomials of degree $d$. 
 Taking the
projectivization, and pulling back to $\p N\vez
W$, we get
$$
\xymatrix
{&\cl O_{\p{N}_{\vphantom{|}}}(-1)
\,\ar@{>->}[d]\ \ar@{->}[dr]^{\ov\rho}& 
\\\ \ 
\cl D_d\ \ \ar@{>->}[r]&
\s d\ \ar@{->>}[r]&
\cl E_d
}
$$
By construction, \ov\rho\ vanishes precisely
over \wtn. This shows that we actually get
\be\label{proj}
\wtn=\p{}(\cl D_d).
\ee
Since the rank of $\cl E_d$ and codimension of \wtn\ 
agree, it follows that \wtn\ represents the top
Chern class of $\cl O_{\p N}(1)\otimes\cl E_d$
(cf.\,\cite[3.2.16,\,p.\,61]{GRR}).
This is the key to the calculation of degrees below.
The map 
\be\label{wtnl}\ba{rcl}
W\vez\p N\supset
\wtn&\stackrel{p_2}\lar&\nld\subset\p N\\
(C,F)&\mapsto& F\ea
\ee
is generically injective
by Noether-Lefschetz theory \cite{Lopez},
cf.\,Corollary\,\ref{maincor} below.
Therefore 
the degree of $ \nld$
can be computed upstairs. Namely, setting
\\\centerline{$m=\dim\wtn,
H=c_1\cl O_{\p N}(1)$,}
we have
$$
\deg \nld=\int{}H^m\cap\wtn=
\int{}H^mc_\tau(\cl O_{\p N}(1)\otimes\cl E_d),
$$
where $\tau=$\,rank$\cl E_d$.
Expanding the top Chern class and pushing
forward to $W$, we arrive at
\be\label{degnld}
\deg\nld=\int_Wc_w(\cl E_d),
\ee
with $w=\dim W$. Since $\cl E_d$ is the pushforward
of a sheaf on $W\vez\p3$, we may apply
Grothendieck-Riemann-Roch \cite[p.\,286]{GRR} 
to express the Chern character of $\cl E_d$ as
$$\ch(\cl E_d)=\ch(q_{1!}(\cl O_{\wt C}(d))=
q_{1\star}\left(\ch(\cl O_{\wt C})\ch(\cl O_{\p3}(d))
\rm todd\p3\right).
$$
Notice that the right hand side is a polynomial in
$d$ of degree $\leq3$. Since $c_w$ is a polynomial of
degree $w$ on the coefficients of the Chern
character, we deduce that $c_w(\cl E_d)$ is a
polynomial in $d$ of degree $\leq3w$.
 \qed

\medskip

\subsection{Remarks.}
\label{cmr}
(1)
The assertion that $\cl E_d$,
as defined in \rf{ed}, is a vector bundle of rank
$p_W(d)$ holds for all $d$ beyond $d_0=$ the maximal
Castelnuovo-Mumford regularity of the members of $W$.
For instance, if $W$ is the family of lines in \p3,
then $d_0=1$.  

\noindent(2)
For the case of elliptic quartic curves
presented below, we note that the regularity of the
ideals \id{ x_{1}^{2}, x_{2}^{2} } and 
\id{x_{1}x_{2}, x_{1}^{2}, x_{2}^{3} }  is $3$, whereas
for \id{x_0^2,x_0x_1,x_0x_2^2,x_1^4} it is $4$. 
The last two ideals  are representatives of the closed
orbits in $W$. An argument of semi-continuity shows
that  $d_0=4$ works for all members of $W$,
see\,\cite{av}. Nevertheless, the map in (\ref{wtnl}) is generically injective only for
$d\geq5$, cf. (\ref{maincor}) below.
Notice that the full Hilb$^{4t}\p3$
has a ``ghost'' component with  regularity
index $\geq5$, see\,\cite{gg}.

\medskip

\section{elliptic quartics}

We consider now the case of surfaces of degree
$\geq4$ containing an elliptic quartic curve in \p3. 
Thus, a general member $C_4$ of $W$ is the
intersection of two quadric surfaces.  The parameter
space $W$ is described in \cite{av} and has been used
in \cite{ell&stromme} to enumerate curves in Calabi-Yau
3-folds. For the convenience of the reader, we
summarize below its main features.

The  Noether--Lefschetz locus of quartic surfaces
containing some $C_4$ is 
slightly exceptional. This is a case when
the map (\ref{wtnl}) fails to be generically injective
(cf. Corollary\,\ref{maincor}): it actually shrinks
$\dim\wtn=34$ to $\dim\nld=33$.  Indeed, if a quartic
surface $F$ contains some general elliptic quartic
 $C_4$, then $F$ must
contain the two pencils $|C_4|$ and $|C_4'|$, where
$C_4'$ is the residual intersection of $F$ with a
quadric containing $C_4$, \ie $2H=C_4+C_4'$, with 
$H=$
plane section. We show in \S\ref{d=4} that $NL(W,4)$
is a hypersurface of degree {\bf38475} in
$\p{34}=|\cl O_{\p3}(4)|$. 

\subsection{}
\label{W4ell}
Next we give an outline of the calculation. 
Put 
$$\x=\G(2,\s2),
$$
the Grassmannian of pencils of quadrics in
\p3. 

The diagram below summarizes the
construction of $W$ as explained in \cite{av}.
\be\label{YZ}
\xymatrix
@C1pc
{
\G(19,\s4)
&\hskip-4cm\supset\hskip-4cm&
W\,\ar@{=}[r]&\wh\x_{\vphantom|}\ \ar@{->}[d]
  &&&&\ \wh\E_{\vphantom|} \ar@{->}[d]\ \ar@{_(->}[llll]
\\
\G(8,\s3)
\vez\x
&&\hskip-4cm\supset\hskip-4cm
&\,\wt\x\ \ar@{->}[d]
  &\hskip-4cm\supset\hskip-4cm&\wt \E
  \ar@{->}[d]&\hskip-4cm\supset\hskip-4cm&\ \,\wt Y\
\ \ar@{->}[d] 
\\
&&\G(2,\s2)\ \ar@{=}[r]
&\x  &\hskip-1cm\supset\hskip-1cm&
Z&\hskip-4cm\supset\hskip-4cm& Y 
}\ee
where 
$$\left\{\ba l
Z\isom\pd3\vez\G(2,\s1)\ \text{ consists of
pencils with a fixed plane};\\
Y\isom\{(p,\ell)\st{} p\supset\ell\}=\ \text{
  closed orbit of }\,Z;\\
\wt Y\lar Y=\p2- \text{bundle of degree 2
  divisors on the varying }\ \ell\subset p;\\
\wh X =\text{ the blow-up of \wt X \ along 
  \wt Y\ and}\\
\wt X =\text{ the blow-up of $X$  along }
Z.
\ea\right.$$

\centerline{$
Z=\left\{
\escala{1.5}\begin{picture}(1.7,.67)(-.1,0.21)
\put(0,0){\line(1,0){1}
\line(1,1){.5}} 
\put(.5,.5){\line(1,0){1}}
\put(0,0){\line(1,1){.5}}
\put(.75,.205){\line(0,1){.651}}
 \end{picture}
\right\}
\supset
Y=
\left\{
\escala{1.5}\begin{picture}(1.7,.67)(-.1,0.21)
\put(0,0){\line(1,0){1}
\line(1,1){.5}} 
\put(.5,.5){\line(1,0){1}}
\put(0,0){\line(1,1){.5}}
\put(.5,0){\line(1,1){.5}}
 \end{picture}
\right\}
\longleftarrow
\wt Y=
\left\{
\escala{1.5}\begin{picture}(1.7,.67)(-.1,0.21)
\put(0,0){\line(1,0){1}
\line(1,1){.5}} 
\put(.5,.5){\line(1,0){1}}
\put(0,0){\line(1,1){.5}}
\put(.5,0){\line(1,1){.5}}
\put(.875,0.44){$\star$}
\put(.374,-.08){$\star$}
 \end{picture}
\right\}
$}

\medskip
Let 
$$
\cl A\subset\s2\vez\x
$$ 
be the tautological subbundle of rank 2
over our Grassmannian of pencils of quadrics.
There is a natural map of vector bundles over
\x\ induced by multiplication,
$$
\mu_3:\cl A\otimes\s1\lar\s3\vez\x
$$
with generic rank 8. It drops rank precisely over
$Z$. It induces a rational map
$\kappa:\x\rar\G(8,\s3)$.  Blowing up \x\ along $Z$, we
find the closure $\wt\x\subset\G(8,\s3)\vez\x$ of the
graph of $\kappa$. Similarly, up on \wt\x\ we have a
subbundle $\cl B\subset\s3\vez\wt\x$ of rank 8 and a
multiplication map
$$
\mu_4:\cl B\otimes\s1\lar\s4\vez\wt\x
$$
with generic rank 19.  The scheme of zeros of
$\stackrel{19}{\bigwedge}\mu_4$
is equal to \wt Y. Indeed, it can be checked
that each fiber of \cl B\
is a linear system of cubics 
which 
\bi
\item
either has base locus equal to a curve with Hilbert
polynomial $p_W(t)=4t$  

\item
or is of the form $p\cdot \s2^{\star\,\star}$,
meaning the linear system with fixed component
a plane $p$, and $\s2^{\star\,\star}$
denoting an 8-dimensional space of quadrics
which define a subscheme  of $p$ of dimension 0 and
degree 2. 
\ei

The exceptional divisor $\wh\E$ is a \p8-bundle over
\wt Y. The fiber of $\wh\E$ over $(p,y_1+y_2)\in\wt Y$
is the system of quartic curves in the plane $p$ which
are singular at the doublet $y_1+y_2$. Precisely, if
$x_0,\dots,x_3$ denote homogeneous coordinates on \p3,
assuming $p:=x_0,\ell=\id{x_0,x_1}$, a typical doublet
has homogeneous ideal of the form
\id{x_0,x_1,f(x_2,x_3)}, with $\deg{} f=2$. Our system
of plane quartics lies in the ideal $\id{x_1,f}^2=
\id{x_1^2,x_1f,f^2}$. Given a non-zero quartic $g$ in
this ideal, we may form the ideal
$J=\id{x_0^2,x_0x_1,x_0f,g}$. It can be checked that
 $J$ contains precisely 19 independent quartics and
its Hilbert polynomial is correct.  In fact, any such
ideal is 4-regular (in the sense of
Castelnuovo-Mumford).  Moreover, up on \wh\x\ we get a
subbundle 
$$\cl
C\subset\s4\vez\wh\x
$$ 
of rank 19. Each of its fibers over \wh\x\ is a system
of quartics that cut out a curve with the correct
Hilbert polynomial.  The multiplication map
$$
\cl C\otimes\s{d-4}\lar\s{d}\vez\wh\x
$$
is of constant rank $\binom{d+3}3-4d$. The image 
$$
\cl D_d\subset\s{d}\vez\wh\x 
$$
is a subbundle as in (\ref{dd}). We have 
$$
\wtn=\p{}(\cl D_d)
\subset\p N\vez\wh\x.
$$
Now the map $\wtn\ra\nld$ is generically injective
for $d\geq5$ in view of Corollary\,\ref{maincor}\,(ii)
below.
Hence the degree of the image
$\nld\subset\p N$ is given by  $\int
c_{16}\cl E_d$, cf.\,(\ref{degnld}).

The above description suffices to feed in
Bott's localization formula with all required data. Indeed,
\wh\x\ inherits a \C\us-action, with 
(a lot) of isolated fixed points. 
The vector bundle $\cl D_d\ra\wh\x$ is
equivariant; ditto for $\cl E_d$. Bott's formula
reads \cite{eg},

\be\label{B}\int
c_{16}\cl E_d=
\sum_{p\in\rm fixpts}
\frac{c_{16}^T( \cl E_d)_p}{c_{16}^T\Tau_p}\,\cdot
\ee

The equivariant classes on the r.h.s are calculated
in two steps. Below we set $d=5$ for simplicity.
First, find the \C\us--weight decomposition of the
fibers of the vector bundles $\cl E_d$ and $\Tau$ at each
fixed point. Then, since the fixed points are isolated,
the equivariant Chern classes $c_i^T$ are just 
the symmetric functions of the weights.

For instance, for the tangent bundle \Tau, say at 
the fixed point corresponding to the pencil 
$p=\id{x_{{0}}^2,x_{{1}}^2}\in\x$, we have
\\\centerline{
$\Tau_p=\Hom(p,\s2/p)=p\ve\otimes\s2/p
=\frac{x_0x_1}{x_0^2}
+\frac{x_0x_2}{x_0^2}+\frac{
x_1x_2}{x_0^2}+
\cdots+\frac{x_3^2}{x_1^2}\,\cdot
$}
Each of the 16 fractions,
 $\frac{x_ix_j}{x_kx_l}$, on the right hand
side symbolizes a 
1-dimensional subspace of $\Tau_p$ where \C\us \  acts
with 
character $t^{x_i+x_j-x_k-x_l}$. The
denominator in \rf{B} is the corresponding
equivariant top Chern class, to wit
\\\centerline{$ 
( x_1-x_0 )  (x_2 -x_0 ) 
 ( x_1+x_2-2x_0)  
\cdots ( 2x_3-2x_1).
$}
The numerator in \rf{B} requires finding  the monomials of
degree 5 that survive modulo the ideal
\\\centerline{$\id{x_0^2,x_1^2}\cdot\id{x_0,\dots,x_3}^3
.$ }
 We are left with 20 (=rank\,$\cl E_5$)
terms,
\\\centerline{$
\ba c
x_2^3x_1x_0+x_3x_2^3x_1+
\cdots+x_3^4x_0.
\ea$}

Now the  equivariant Chern class $c_{16}^T(\cl E_d)_p$
is the coefficient of $t^{16}$ in the product 
\\\centerline{$(t+x_0+x_1+3x_2)(t+x_1+3x_2+x_3)
\cdots(t+x_0+4x_3)$}.
(20 factors.)
In practice, all these calculations are made
substituting $x_i$ for suitable numerical
values, cf. the computer algebra  scripts in
\cite{v}.

\subsection{the case $d=4$}\label{d=4}
Presently $p_2:\wt{NL}(W,4)\ra{NL}(W,4)$
is no longer generically injective. It shrinks
dimension by one: a general fiber is a disjoint
union of \p1's (cf. Corollary\,\ref{maincor}\,(iii) below).
Explicitly, say $F_4=A_1Q_1+A_2Q_2,\,\deg
A_i=\deg Q_i=2$, with $A_i$ and $Q_i$ generic.
Then 
\\\centerline{$
F_4=(A_1-tQ_2)Q_1+(A_2+tQ_1)Q_2$,}
so $F_4$ contains the pencil of elliptic quartics
$
\id{A_1-tQ_2,A_2+tQ_1},t\in\p1$;
setting $t=\infty$, we find \id{Q_1,Q_2}. Similarly, get 
\id{Q_1-tA_2,Q_2+tA_1}. This is one and the
same pencil. But there is also \id{A_1-tA_2,Q_2+tQ_1}.
In general, these 2 pencils are disjoint. Looking
at them as curves in $\x=\G(2,\s2)$, we actually
get a Pl\"ucker--embedded conic,
$(A_1-tA_2)\wedge(Q_2+tQ_1)= 
A_1\wedge{}Q_2+t(A_1\wedge{}Q_1-A_2\wedge{}Q_2)
-t^2A_2\wedge Q_1$,
disjoint from $Y$ (see \rf{YZ}). In
particular, capping each conic against the Pl\"ucker
hyperplane class
$\Pi=-c_1\cl A$, we find 2. As before, we may write
$$
\deg \nld=\int{}H^{33}\cap\nld.
$$ 
The cycle $p_2\us H^{33}\cap\nld$ can be
represented by a sum of $\deg \nld$ disjoint
unions of pairs of \p1's.
Hence 
$$
\ba r
\deg \nld
=\ds{\frac14}\ds{\int_{\wh\x}}\Pi{\cdot}H^{33}\cap\wtn
=\frac14\int_{\wh\x}\Pi{\cdot}(p_1)\ls{}H^{33}\cap\wtn
\\\na9
=\ds{\frac14\int_{\wh\x}}\Pi{\cdot}c_{15}(\cl E_4).
\ea$$
The latter integral can be computed via Bott's
formula and we get {\bf38475}, cf.\, the script in \cite{v}.
This has been found independently in \cite{kr} with
different techniques, using \cite{rp}.

\section{The fibers of $p_2$}
The main result needed to validate the above enumeration
is  the following.

\subsection{Proposition}
{\em Let $C \subset \p3$ be a smooth irreducible curve
of degree $e$ and genus $g$. Let $d \gg 0$ and let $F
\subset \p3$ be a general surface of degree $d$
containing $C$. Then $C$ is the only  effective
divisor of degree $e$ and arithmetic genus $g$ on $F$. 
}

\begin{proof}
By \cite[Cor.II.3.8]{Lopez} we have that $\Pic(F)$ is
freely generated by its hyperplane section $H$ and
$C$. Let $C'$ be an effective divisor of degree $e$ and
arithmetic genus $g$ on $F$. Then there are two
integers $a, b$ such that, on $F$, we have $C' \sim aH
+ bC$. Now $e = H\cdot{}C' = ad + b e$, so that
\begin{equation}
\label{uno}
a = \frac{e}{d}(1-b).
\end{equation}
By adjunction formula we have 
\be\label{c2}
C^2 = 2p_a(C) - 2 - K_F\cdot{}C = 2g-2-e(d-4)
  = (C')^2 < 0
\ee
 as $d \gg 0$.
Now
\[ C^2 = (C')^2 = (aH + bC)^2 = a^2d + b^2 C^2 + 2abe \]
and using \eqref{uno} we get
\[ (1-b^2) (e^2 - dC^2) = 0.
\]
Note that $e^2 - dC^2 > 0$ by (\ref{c2}), 
whence $b = \pm 1$. If $b =
-1$ we have from (\ref{uno}) 
$a = \frac{2e}{d} \not\in \bbb Z$, as $d
\gg 0$. Therefore we deduce that $b = 1$ and $a = 0$,
that is $C' \sim C$. Since $C^2 < 0$ we must have $C' =
C$.
\end{proof}

\subsection{Corollary}\label{maincor}
{\em 
Let W be an irreducible subvariety of a  Hilbert scheme
component of  curves in $\p3$  of degree e and
arithmetic genus g with general member smooth. 
Let $\p N = \p{}(
H^0(\cl O_{\p3}(d)))$. Then 

{\rm(i)}
there is a $d_0$ such that
for all  $d\geq d_0$
the projection map
 $p_2 : \wt{NL}(W, d) \to \p N$  is generically
 injective.

{\rm(ii)} If 
W is the family of elliptic quartics then we can take 
$d_0= 5$, \ie $p_2$ is generically one-to-one for 
$d\geq5$ and

{\rm(iii)} 
for $d = 4$, the general fiber of $\wt{NL}(W,4)
\stackrel{p_2}\lar{NL}(W,4)$ is two disjoint
$\p1$'s.}

\begin{proof}

We know from (\ref{proj}) that $p_1:\wtn\ra W$ is a
projective bundle. Hence
$NL(W,d) = p_2(\widetilde{NL}(W, d))$ is irreducible
and a general element $F \in NL(W, d)$ can be
identified with a general hypersurface of
degree $d$ containing a general $C\in W$. Hence
assertion  (i) follows from the propostion.
Assertion (ii) also follows, except for $d=8$. In this
case, 
with the notation as in (\ref{uno}), if $b=-1$
we would get $a=1$, whence $C' \sim H
-C$ so that $C$ would be contained
in a plane, which is absurd. If $b=1$, may proceed as at the end of
the proof of the above proposition.

If $d = 4$ then we get instead $C' \sim 2H - C$.
The exact sequence
\[ 0 \lar \cl O_F \lar \cl O_F(C) \lar \cl O_C (C)
=\cl O_C \lar
0 \] shows that $|\cl O_F(C)| \cong \p1$ and
similarly $|\cl O_F(2H-C)| \cong |\cl O_F(C')| \cong
\p1$. Moreover there is no curve $D$ on
$F$ such that $D \sim C$ and $D \sim
2H - C$ for then $2C \sim 2H$ giving the
contradiction $0 = C^2 = H^2 = 4$. This proves
that, in this case, the general fiber of $p_2$
is two disjoint $\p1$'s,
\begin{equation}
\label{quattro}
p_2^{-1}(F) = |\cl O_F(C)| \cup |\cl O_F(2H-C)|.
\end{equation}
\end{proof}

\subsection{The formula}\label{f}
In view of Proposition\,\ref{pol}, it suffices to
find the degrees of \nld\ for $3\cdot16$+1 values of
$d\geq5$ and interpolate. This is done in
\cite{v}. We obtain

{\small
$$\ba c
\binom{d-2}3
\Big(106984881d^{29}-3409514775d^{28}+
57226549167d^{27}-
\\643910429259d^{26}+5267988084411d^{25}-31628193518727d^{24}+
126939490699539d^{23}
\\-144650681793207d^{22}-2701978741671631d^{21}+
28913126128882647d^{20}-
\\182919422241175163d^{19}
+858473373993063183d^{18}-
3061191057059772423d^{17}
\\+7448109470245631187d^{16}-3841505361473930575d^{15}-
80644842327962348733d^{14}+
\\568059231910087276234d^{13}-2560865812030993315212d^{
12}+
9159430737614259196104d^{11}
\\-27608527286339077691280d^{10}+
71605637662357479581024d^{9}
\\-160009170853633152594240d^8
+303685692157317249665152d^7
\\-473993548940769326728704d^6+
571505502502703378479104d^5
\\-459462480152611231457280d^4+
111908571251948243582976d^3
\\+251116612534424272896000d^2-
328452832055501940326400d
\\+136886449647246114816000\Big)
\Big/
(2^{27}\cdot3^{9}\cdot5^2\cdot7^2\cdot11\cdot13).
\ea$$
}

\subsection{Acknowledgments}
The authors thank the referee for kindly pointing out
to them a sizable portion of the argument
that required clarification. They are also grateful to
Kristian Ranestad and 
Rahul Pandharipande for correcting an error in 
\S\ref{d=4} and for calling their attention to 
\cite{kr} and \cite{rp}.

\section*{appendix}
Script for Singular.
This builds on the original script for Maple found in
P. Meurer, \cite{Meurer}. Copy and paste it into a
Singular session.

\begin{verbatim}
//<"frac4bott";
option (noloadLib,noredefine);
LIB "ring.lib";

proc Proj(n){
execute("ring P"+string(n)+"=0,x(0..n),dp;
def xx=ideal(x(0..n));export xx;
keepring P"+string(n)); nameof(basering);basering;}

proc omit(w,a)
{list ll,lw,v;int i,j,s,t,ch;
 s, t=size( w),size( a);
    if(typeof( a)<>"list") {
list ll;
if(typeof(a)=="poly") {
t=1;ll[1]=a;}
if(typeof(a)=="ideal") {
t=size(matrix(a));ch=1;}
if(typeof(a)=="intvec" or typeof(a)=="int")
{t=size((a));ch=1;}
if(ch==1){
for( i=1; i<= t; i++){ll[i]=a[i];}}}
else{ll= a;}
    if(typeof(w)<>"list"){
list lw;
if(typeof(w)=="poly") {
s=1;lw[1]=w;}
if(typeof(w)=="ideal") {
s=size(matrix(w));ch=1;}
if(typeof(w)=="intvec" or typeof(w)=="int")
{s=size((w));ch=1;}
if(ch==1){
for( i=1; i<= s; i++){lw[i]=w[i];}}
def ty=typeof(w);
}else{lw=w;}
    for(i=1;i<=s;i++)
    {ch=0;
for(j=1;j<=t;j++){
	if(lw[i]==ll[j]){ch=1;break;}}
        if(ch==0){
v=insert(v,lw[i],size(v));}
    }if(defined(ty)<>0){
	if(size(v)<>0){def st=string(v);}
	else{def st="0";}
	execute(ty+" vv="+st);}
	else{def vv=v;}
return(vv);}

proc origin(v,p)
{def q=p;if(typeof(q)<>"list"){
 for(int i=1;i<=size(v);i++){q=subst(q,v[i],poly(0));}
 return(q);}
else{for(int i=1;i<=size(q);i++){q[i]=origin(v,q[i]);}
return(q); } }

proc hilbp//(ideal I)
{int nargs=size(#);ideal I=#[1];
  //start comparing rings
  poly pvars=indets(I);
  if(pvars==0){pvars=indets(vars(1..(nvars(basering))));}
  if(pvars<>indets(vars(1..(nvars(basering)))))
    {if(nargs==3)
      {print("too many vars in basering "+nameof(basering)+
	     ", will change to less");}
    if (nargs>=2){pvars=indets(#[2]);}
string st="ring r0=(0";if(npars(basering)<>0){st=st+","+parstr(basering);}
	st=st+"),(";
    for(int i=1;i<=(size(pvars)-1);i++){st=st+string(pvars[i])+",";}
st=st+string(pvars[i])+"),(dp)";
 st=st+";ideal I="+string(I);
 execute(st);if(nargs==3){where();}}
  intvec iv=hilbPoly(I);if(nargs==3){print(iv);}
 int s=size(iv);string st,st1;
 for(int i=1;i<=s;i++){
   if(iv[i]<>0){st1=string(iv[i]);if(i==1){st=st1;}
   else{if(i==2){if(iv[i]<>1){st1=st1+"*t";}
   else{st1="t";}if(size(st)>0){if(iv[i]<0){st=st+st1;}else{st=st+"+"+st1;}}
   else{st=st1;}}
   else{if(iv[i]<>1){st1=st1+"*t^"+string(i-1);}
   else{st1="t^"+string(i-1);}if(size(st)>0){
     if(iv[i]<0){st=st+st1;}else{st=st+"+"+st1;}}
   else{st=st1;}}}}}
 if(s>=3){st="(" + st + ")/" + string(factorial(s-1));}
 return(st); }

proc indets(p)//returns sum of variables
{int i,j,s;poly v,w,v1;ideal pp;
if (typeof(p)=="int" or typeof(p)=="intvec"
 or typeof(p)=="intmat"  or typeof(p)=="number")
{v=0;return(v);}
else{
if(typeof(p)=="poly")   {
	for (i=1;i<=nvars(basering);i++)
          {if(reduce(p,std(var(i)))<>p)
            {v=v+var(i);}
          }
return(v);}
else{if(typeof(p)=="ideal") {
   pp=simplify(p,2);s=size(pp);w=0;
   if(s==0){v=0;return(v);}
   else{
     for(i=1;i<=s;i++)
          {v=indets(pp[i]); 
      if(v<>0){v1=leadmonom(v);
         w=w+v*v1^i;}};
return(indets(w));}}
else
{if (typeof(p)=="vector" or typeof(p)=="matrix")
  {return(indets(ideal(p)));
  }
else
{if (typeof(p)=="list")
  {def s=size(p);if (s==0){return(poly(0));}
else{ideal pp;
for(i=1;i<=s;i++){pp[i]=indets(p[i]);
}return(indets(pp));}}
}}}}}

proc vars(v)
{ideal z;def vv=v;if (typeof(vv)=="intvec") {
for (int i=1;i<=size(vv);i++){z=z+var(vv[i]);}
 return(z);}
 else{if (typeof(vv)=="ring") {
  if(string(vv)<>string(basering)) {
  def goback=nameof(basering);
  setring vv;def p=vars(intvec(1..nvars(vv)));
  execute("setring "+goback);def p=imap(vv,p);
 return(p);}else{return(vars(intvec(1..nvars(vv))));}
}}}

proc mdc
{int s,i;s=size(#);string ty=typeof(#[1]);
 if(s==1){
        if (ty=="poly" or ty=="int" ){return(#[1]);}
        else{
                if(ty=="list" or ty=="vector" or
                   ty=="ideal" or ty=="intvec"
           or ty=="matrix" )
                        {if(ty=="vector" or ty=="matrix")
                          {ideal ii=ideal(#[1]);
                          }
                         else{
                                execute(ty+" ii=#[1];");
                             }
                        if (ty<>"intvec"){
                        ii=simplify(ii,2);}
                        i=size(ii);
            if(i==0){return(poly(0));}
            else{return(mdc(ii[1..i]));
            }           
                          }
                }
            }
 if (s==2){
                return(gcd(#[1],#[2]));
          }
 if (s>2){
        execute(ty+" p=#[1];");
        p=gcd(#[1],#[2]);
        for(i=3;i<=s;i++){p=gcd(p,#[i]);}
return(p);}}

proc dotprod( w1, w2) {
   if(typeof(w1)=="list"
      & typeof(w1)=="list")
     {int s=min(size(w1),size(w2));poly p;
  for (int j=1;j<=s;j++){p=p+w1[j]*w2[j];}
  return(p);}
else{list v1,v2;
if(typeof(w1)=="vector") {//"//w1 vector";
for(int i=1;i<=nrows(w1);i++){v1[i]=w1[i];}}
 else{ if(typeof(w1)=="ideal") {//"//w1 ideal";
for(int i=1;i<=ncols(w1);i++){v1[i]=w1[i];}}
   else{if(typeof(w1)=="intvec"){//"//w1 intvec";
for(int i=1;i<=size(w1);i++){v1[i]=w1[i];}}
     else{ERROR(typeof(w1)+
"?w1?not (int)vector nor ideal");}}}
if(typeof(w2)=="vector") {//"//w2 vector";
for(int i=1;i<=nrows(w2);i++){v2[i]=w2[i];}}
 else{ if(typeof(w2)=="ideal") {//"//w2 ideal";
for(int i=1;i<=ncols(w2);i++){v2[i]=w2[i];}}
   else{if(typeof(w2)=="intvec"){//"//w2 intvec";
for(int i=1;i<=size(w2);i++){v2[i]=w2[i];}}
     else{ERROR(typeof(w2)+
"?w2?not (int)vector nor ideal");}}}
 return(dotprod(v1,v2));}}

proc min//intvec or list or whatever?
{def s=size(#);
 def m=#[1];for(int i=2;i<=s;i++){if(#[i]<m){m=#[i];}}
return(m); }

proc prod//(ideal y)
{ def s,ty=size(#),typeof(#[1]);
  poly p=1;int i;
  if(ty=="ideal") {def y=#[1];
  matrix m=matrix(y);//bypass size problem
  for(i=1;i<=ncols(m);i++){p=p*m[1,i];}
  }
  else{
    for(i=1;i<=s;i++){p=p*#[i];}}
  return(p);}

proc coefmon(poly p, poly m)
//coeff of monom m in poly p
{ if(m==1 or p==0){return(jet(p,0));}
 else{ideal m1=std(m);def v=std(pol2id(indets(m)));poly q0;
 for(int i=1;i<=size(p);i++){def q=reduce(p[i],m1);
 if(q==0){q=p[i]/m;
 if(reduce(q,v)==q) {q0=q0+q;}}}
   return(q0);} }

proc pol2id(poly p){
ideal I;for(int i=1;i<=size(p);i++)
        {I[i]=p[i];} return(I);}

proc sumfrac//(list ab,list cd,...)
{def nargs=size(#);//"//"+string(nargs);
if (nargs==2) {list ab=#[1];list cd=#[2];
def l= (ideal(ab[1]*cd[2]+ab[2]*cd[1],
ab[2]*cd[2]));def p=gcd(l[1],l[2]);
return(list(l[1]/p,l[2]/p));}
else{def ll=sumfrac(#[1],#[2]);
for(int i=3;i<=nargs;i++){
ll=sumfrac(ll,#[i]);} return(ll); }}

proc subfrac(list ab,list cd)
{def l= (ideal(ab[1]*cd[2]-ab[2]*cd[1],
ab[2]*cd[2]));def p=gcd(l[1],l[2]);
return(list(l[1]/p,l[2]/p));}

proc mulfrac(list ab,list cd)
{def l=  ideal(ab[1]*cd[1],ab[2]*cd[2]);
def p=gcd(l[1],l[2]);
return(list(l[1]/p,l[2]/p));}

proc invfrac(list ab)
{return(list(ab[2],ab[1]));}

proc divfrac(list ab,list cd)
{return(mulfrac(ab,invfrac(cd)));}

proc mon(vs,d)
 {if(typeof(vs)=="poly"){def vvs=pol2id(vs);}
else{def vvs=vs;}
return(vvs^d);}  

proc mylead(poly h){
def oldring=nameof(basering);
def goback="setring "+oldring;
def h1=h;
if(find(string(h1),"x")<>0){setring rx;}
else{setring rz;def xx=zs;}
execute("def f=imap("+oldring+",h1)");
def iv,a=leadexp(f),leadcoef(f);
a=int(a);
f=(dotprod(iv,xx));
def a1,h1=string(a),string(f);
if(find(h1,"x")==0) { setring rxz;
execute("def f="+h1);
f=z2x(f);h1=string(f);setring rx;
execute("def f="+h1);
}
f=x2p(f);list mm;
execute("int j="+a1);
for(int i=1;i<=(j);i++){mm[i]=f;}
if(oldring<>"rx"){execute(goback);
def mm=imap(rx,mm);}
return(mm); }

proc pesos(H){
//H= poly or list(num,den) of fraction
list rres;
def ty=typeof(H);
if(ty=="poly" or ty=="number" or ty=="int"){
if(indets(H)==0){return(list(0));}else{
rres=mylead(H[1]);
   for(int i=2;i<=size(H);i++){
rres=rres+mylead(H[i]) ;}
return(rres);}}
else{if(ty=="list" and size(H)==2) {
def num,den=H[1..2];
//den assumed monomial
rres=pesos(H[1]);def denn=pesos(den);
for(int i=1;i<=size(rres);i++) {
rres[i]=rres[i]-denn[1];}
return(rres);}
else{if(ty=="ideal")
{rres=pesos(H[1]);
for(int i=2;i<=size(matrix(H));i++) {
rres=rres+pesos(H[i]);}} return(rres);}}}

proc cherns(H){
if(defined(t)<>-1){def t=var(1);}
poly sigma= 1;
for(int j=1;j<=size(H);j++){
sigma=reduce((1+H[j]*t)*sigma,std(t^(DIM+1)));}
list l;
for (int j=1;j<=DIM;j++){
l[j]=coefmon(sigma,t^j);}
return(l);}

proc topchern(H) {return(prod(H));}

proc TGrass(idd_, E_){
if(typeof(idd_)<>"list") {
def En_,iddn_=list(E_,1),list(idd_,1);}
else{def iddn_,En_= idd_,E_;}
def p=subfrac(En_,iddn_);
def L=divfrac(p,list(iddn_[1][1],iddn_[2]));

for(int j=2;j<=size(iddn_[1]);j++){
L=sumfrac(L,divfrac(p,list(iddn_[1][j],iddn_[2])));}
return(L);}

proc x2p(f){
  if(defined(n)==0){return("//need int n");}
  if(defined(ivP)==0){
    return("//need intvec ivP=0,2,7,10...");}
  if(defined(x(n))==0){return(
"//need  Proj(n); rx=Pn;ideal xx=x(0..n)");}
def g=subst(f,x(0),poly(ivP[1]));
for(int i=1;i<=n;i++){
g=subst(g,x(i),poly(ivP[i+1]));}
return(g); }
//fixpts start here
//<"4ellsing"; <"frac4bottarXiv";
int n=3; Proj(n);def rx=P3;ring r=0,t,dp;r=r+rx;setring r;
imapall(P3);
def xx2=xx^2;def xx2s=sum(xx2);
def xx4=std(xx2^2);def xx4s=sum(xx4);
def xs=sum(xx);

intvec ivP=0, 1, 5, 18;int DIM=16;
list G;

for(int i=1;i<=size(xx2);i++){
for(int j=i+1;j<=size(xx2);j++){
G[1+size(G)]=list(list(xx2[i],xx2[j],
"G:q1,q2,tg"),
TGrass(xx2[i]+xx2[j],xx2s)
);}}

list Z,G2; //Z=P3* x G(2,4), plane, line
for(int i=1;i<= size(G);i++){
def l=G[i]; def q1,q2=l[1][1],l[1][2]; 
def p=gcd(q1,q2);
if( p<>1 ) {def l1,l2=q1/p,q2/p; def 
tg=sumfrac(TGrass(l1+l2,xs),TGrass(p,xs));
def nml=subfrac(l[2],tg);
Z[1+size(Z)]
=list(list(p,l1,l2,i,"Z:p=mdc,l1,l2,pos.G,tg,nml"),
tg,nml);}
else{def q4=ideal(q1,q2)*xx2;
G2[1+size(G2)]
=list(list(q1,q2),G[i][2],q4)  ;
if( size(q4)<>19 ){ string(i);break;}}}

list E1,G2E1,W; //exc over Z
//should get welldef system of 8 cubics everywhere
//W=P2-bdle over clsd orbit of Z
for(int i=1;i<= size(Z);i++){
  def l=Z[i]; def L1=G[l[1][4]][1];L1=L1[1..2]; 
  nml=l[3];
  def num,den=nml[1..2];
    for(int i1=1;i1<= size(num);i1++) {
    def exc=sumfrac(list(0,1),list(num[i1],den));
    int j1=0;
for(int j=1;j<=size(L1);j++){ def ll=L1[j];
def j0=ll/exc[2];
if(j0<>0 ){j1=1;
//check div instead
break; }}
if(j1==1){
def L2=omit(L1,ll);
L2=L2[1],ll+t*exc[1]*j0;
def qs=sat(ideal(L2[1..2])*xx,t)[1];
qs=reduce(qs,xx4);
qs=std(origin(t,qs));def j0=size(qs),hilbp(qs);
    tg=divfrac(subfrac(nml,exc),exc);//(nml-exc)/exc
    tg=sumfrac(tg,l[2]);//tg fiber+base
    tg=sumfrac(tg,exc);def ll=l[1];
    E1[1+size(E1)] =list(list(ll[1..4],i,
"E1: mdc,l1,l2,pos.G,pos.Z,2:tg,3:exc,4:8 cubics"),
tg, exc,qs); 
if (j0[2]<>"4*t") { def L=mdc(qs); //L=l[1][1]
def qs0=ideal(qs/L); def l1=l[1][2]+l[1][3]-L;
def tg0=sumfrac(TGrass(L,xs) //tg P3dual
,TGrass(l1,xs-L));  //tg fiber W0-->P3dual
def q0=origin(ideal(L,l1),list(sum(qs0),xx2s));
tg0=sumfrac(tg0,TGrass(q0[1],q0[2]));
    W[1+size(W)]=list(list(L,l1,i,size(E1),
"W:mdc,l1,pos.Z,pos.E1,2:tgW,3:nml,4:exc,5:8th q,6:8 qs"),
tg0, subfrac(tg,tg0), exc,q0[1],qs0); }
else {def qs4=std(qs*xx);
G2E1[1+size(G2E1)]
=list(list(size(E1),tg,qs4));
if( size(qs4)<>19 ){ "//?"+string(i);}
}}
else//( j1==0 )
{ "//??"+string(i,",",i1);}//never happens, no need 
}}

def h = "4*t";
list E2; //exc over W; //P8-bdle over W; dimW=3+2+2
for(int i=1;i<= size(W);i++) {
def  l=W[i]; l=l[1];
def L,l1,q0=l[1..2],W[i][5];
def nml=W[i][3];
  def num,den=nml[1..2];
def qs=E1[l[4]][4];
qs=std(qs*xx);
  for  (int i1=1;i1<= size(num);i1++) {
    exc=sumfrac(list(0,1),list(num[i1],den));
    tg=divfrac(subfrac(nml,exc),exc);//(nml-exc)/exc
    tg=sumfrac(tg,W[i][2]);//tg fibra+base
    tg=sumfrac(tg,exc);
    def qs1= qs + num[i1];
E2[1+size(E2)]=list(list(L,l1,q0,i,
"E2: L,l1,q0,pos.W,2:tg,3:exc,4:19 qtics"),
		tg, exc, qs1);
}}

def l=list(G,Z,E1,E2,G2,G2E1);
setring rx;
def l=imap(r,l);
def G,Z,E1,E2,G2,G2E1=l[1..size(l)];

//actual numerical computations at fixpts start here

//case d=4

poly f;def X=G2;
for(int i=1;i<=size(X);i++) {
def H=-sum(pesos(sum(X[i][1])));//1st chern class
//=Plucker hypln class
def qs,c=std(X[i][3]),topchern(pesos(X[i][2]));
def d=cherns(pesos(kbase(qs,4)))[DIM-1]*H;
f=f+d/c;
if(i mod 10 ==1) {"//in G2, remain "+string(size(X)-i);
system("sh","date");
}}
"//done G2";

X=G2E1;//X[1][1]...
for(int i=1;i<=size(X);i++){
def l=E1[X[i][1][1]][1];
def q12=l[1]*l[2]+l[1]*l[3];
def H=-sum(pesos(q12));//1st chern class
def qs,c=std(X[i][1][3]),topchern(pesos(X[i][1][2]));
def d=cherns(pesos(kbase(qs,4)))[DIM-1]*H;
f=f+d/c;
if(i mod 10 ==1) {"//in G2E1, remain "+string(size(X)-i);
system("sh","date"); } }

X=E2;
for(int i=1;i<=size(X);i++){
def l=X[i];
def q12=l[1][1]^2+l[1][1]*l[1][2];
def H=-sum(pesos(q12));//1st chern class
def qs,c=std(l[4]),topchern(pesos(l[2]));
def d=cherns(pesos(kbase(qs,4)))[DIM-1]*H;
f=f+d/c;
if(i mod 10 ==1) {"//in E2, remain "+string(size(X)-i);
system("sh","date"); } } f;
f=f/4;
"//case d=4: "+string(f);

//choose maxdeg>=5 ; need 48...
int minideg,maxideg=54,54;

def X=G2;poly f(minideg..maxideg);
for(int i=1;i<=size(X);i++){
def qs,c=std(X[i][3]),topchern(pesos(X[i][2]));
for (int i1=minideg;i1<=maxideg;i1++){
def d=cherns(pesos(kbase(qs,i1)))[DIM];
f(i1)=f(i1)+d/c;}
if(i mod 10 ==1) {"//in G2, remain "+string(size(X)-i);
system("sh","date"); } }"//done G2";

X=G2E1;//X[1][1]...
for(int i=1;i<=size(X);i++){
def qs,c=std(X[i][1][3]),topchern(pesos(X[i][1][2]));
for (int i1=minideg;i1<=maxideg;i1++){
def d=cherns(pesos(kbase(qs,i1)))[DIM];
f(i1)=f(i1)+d/c;}
if(i mod 10 ==1) {"//in G2E1, remain "+string(size(X)-i);
system("sh","date"); } }"//done G2E1";

X=E2;
for(int i=1;i<=size(X);i++){
def qs,c=std(X[i][4]),topchern(pesos(X[i][2]));
for (int i1=minideg;i1<=maxideg;i1++){
def d=cherns(pesos(kbase(qs,i1)))[DIM];
f(i1)=f(i1)+d/c;}
if(i mod 10 ==1) {"//in E2, remain "+string(size(X)-i);
system("sh","date"); } }"//done E2";

for (int i1=minideg;i1<=maxideg;i1++){
string(i1,",",f(i1));}

if(1>2){//check answer
def t=var(1);
def 
rr=(t-2)*(t-3)*(t-4)*(106984881*t^29-3409514775*t^28+
57226549167*t^27-643910429259*t^26+5267988084411*t^25
-31628193518727*t^24+126939490699539*t^23-144650681793207*t^22
-2701978741671631*t^21+28913126128882647*t^20
-182919422241175163*t^19+858473373993063183*t^18-\
3061191057059772423*t^17+7448109470245631187*t^16
-3841505361473930575*t^15-80644842327962348733*t^14+
568059231910087276234*t^13-2560865812030993315212*t^12+
9159430737614259196104*t^11-27608527286339077691280*t^10+
71605637662357479581024*t^9-160009170853633152594240*t^8+
303685692157317249665152*t^7-473993548940769326728704*t^6+
571505502502703378479104*t^5-459462480152611231457280*t^4+
111908571251948243582976*t^3+251116612534424272896000*t^2-\
328452832055501940326400*t+136886449647246114816000)
/(number(2)^28*(3)^10*(5)^2*(7)^2*11*13); 
//subs(t=54,rr);
}1>2

\end{verbatim}

\end{document}